\documentclass[12pt]{amsart}
\usepackage{amscd,amssymb,amsthm}
\usepackage{graphicx}
\usepackage{hyperref}

\newtheorem{theorem}{Theorem}[section]
\newtheorem{remark}{Remark}[section]
\newtheorem{corollary}{Corollary}[section]

\newtheorem{example}{Example}[section]
\newtheorem{proposition}{Proposition}[section]
\numberwithin{equation}{section}

\title[To appear in: \textit{Inequality Theory and Applications},
\textbf{4(}2004)]{Relative Divergence Measures \\ and Information
Inequalities}

\author{Inder Jeet Taneja}
\address{Departamento de Matem\'{a}tica\\ Universidade Federal
de Santa Catarina\\
88.040-900 Florian\'{o}polis, SC, Brazil}

\email{taneja@mtm.ufsc.br}

\urladdr{http://www.mtm.ufsc.br/$\sim $taneja}

\thanks{To appear in: \textit{Inequality Theory and Applications},
\textbf{4(}2004)}



\begin{document}

\begin{abstract}
There are many information and divergence measures exist in the
literature on information theory and statistics. The most famous
among them are Kullback-Leiber's \cite{kul} \textit{relative
information} and Jeffreys \cite{jef} \textit{J-divergence},
\textit{Information radius} or \textit{Jensen difference
divergence measure} due to Sibson \cite{sib}. Burbea and Rao
\cite{bur1, bur2} has also found its applications in the
literature. Taneja \cite{tan2} studied another kind of divergence
measure based on \textit{arithmetic and geometric means}. These
three divergence measures bear a good relationship among each
other. But there are another measures arising due to
\textit{J-divergence}, \textit{JS-divergence} and
\textit{AG-divergence}. These measures we call here
\textit{relative divergence measures} or \textit{non-symmetric
divergence measures}. Here our aim is to obtain bounds on
\textit{symmetric }and \textit{non-symmetric divergence measures}
in terms of \textit{relative information of type s} using
properties of Csisz\'{a}r's \textit{f-divergence}.
\end{abstract}

\maketitle

\section{Introduction}

Let
\[
\Gamma _n = \left\{ {P = (p_1 ,p_2 ,...,p_n )\left| {p_i >
0,\sum\limits_{i = 1}^n {p_i = 1} } \right.} \right\}, \,\, n
\geqslant 2,
\]

\noindent be the set of all complete finite discrete probability
distributions.

The Kullback Leibler's (1951) \textit{relative information} is given by
\begin{equation}
\label{eq1}
K(P\vert \vert Q) = \sum\limits_{i = 1}^n {p_i \ln (\frac{p_i }{q_i })} ,
\end{equation}

\noindent for all $P,Q \in \Gamma _n $.

We observe that the measure (\ref{eq1}) is not symmetric in $P$
and $Q$. Its \textit{symmetric version} famous as
\textit{J-divergence} (Jeffreys \cite{jef}; Kullback and Leiber
\cite{kul}) is given by
\begin{equation}
\label{eq2}
J(P\vert \vert Q) = K(P\vert \vert Q) + K(Q\vert \vert P)
 = \sum\limits_{i = 1}^n {(p_i - q_i )\ln (\frac{p_i }{q_i })} .
\end{equation}

For simplicity, we call here the measure $K(Q\vert \vert P)$ the
\textit{adjoint} of $K(P\vert \vert Q)$ and vice-versa.

Alternatively, the measure $J(P\vert \vert Q)$ can also be written
in the following way:
\begin{equation}
\label{eq3}
J(P\vert \vert Q) = D(P\vert \vert Q) + D(Q\vert \vert P),
\end{equation}

\noindent where
\begin{equation}
\label{eq4}
D(P\vert \vert Q) = \sum\limits_{i = 1}^n {(p_i - q_i )\ln \left( {\frac{p_i
+ q_i }{2q_i }} \right)}
\end{equation}

\noindent and
\begin{equation}
\label{eq5}
D(Q\vert \vert P) = \sum\limits_{i = 1}^n {(q_i - p_i )\ln \left( {\frac{p_i
+ q_i }{2p_i }} \right)} .
\end{equation}

Let us consider the following two measures:
\begin{equation}
\label{eq6}
F(P\vert \vert Q) = K\left( {P\vert \vert \frac{P + Q}{2}} \right) =
\sum\limits_{i = 1}^n {p_i \ln \left( {\frac{2p_i }{p_i + q_i }} \right)}
\end{equation}

\noindent and
\begin{equation}
\label{eq7}
G(P\vert \vert Q) = K\left( {\frac{P + Q}{2}\vert \vert P} \right) =
\sum\limits_{i = 1}^n {\left( {\frac{p_i + q_i }{2}} \right)\ln \left(
{\frac{p_i + q_i }{2p_i }} \right)} .
\end{equation}

The \textit{adjoint} forms of the measures (\ref{eq6}) and
(\ref{eq7}) are given by
\begin{equation}
\label{eq8}
F(Q\vert \vert P) = K\left( {P\vert \vert \frac{P + Q}{2}} \right) =
\sum\limits_{i = 1}^n {q_i \ln \left( {\frac{2q_i }{p_i + q_i }} \right)}
\end{equation}

\noindent and
\begin{equation}
\label{eq9}
G(Q\vert \vert P) = K\left( {\frac{P + Q}{2}\vert \vert Q} \right) =
\sum\limits_{i = 1}^n {\left( {\frac{p_i + q_i }{2}} \right)\ln \left(
{\frac{p_i + q_i }{2q_i }} \right)} .
\end{equation}

\noindent respectively, and the \textit{symmetric forms} are given
by
\begin{equation}
\label{eq10}
I(P\vert \vert Q) = \frac{1}{2}\left[ {F(P\vert \vert Q) + F(Q\vert \vert
P)} \right]
\end{equation}

\noindent and
\begin{equation}
\label{eq11}
T(P\vert \vert Q) = \frac{1}{2}\left[ {G(P\vert \vert Q) + G(Q\vert \vert
P)} \right],
\end{equation}

\noindent respectively.

The three measures $J(P\vert \vert Q)$, $I(P\vert \vert Q)$ and
$T(P\vert \vert Q)$ are related with each other as
\begin{equation}
\label{eq12} J(P\vert \vert Q) = 4\left[ {I(P\vert \vert Q) +
T(P\vert \vert Q)} \right].
\end{equation}

Moreover, the measures (\ref{eq4}) can also be written as
\begin{equation}
\label{eq13}
D(P\vert \vert Q) = 2\left[ {F(Q\vert \vert P) + G(Q\vert \vert P)} \right]
\end{equation}

The measure $I(P\vert \vert Q)$ is famous in the literature as
\textit{information radius} (Sibson \cite{sib}) or \textit{Jensen
difference divergence measure}, or simply, \textit{JS-divergence}
(Burbea and Rao \cite{bur1, bur2}). The measure (\ref{eq9}) is new
in the literature and is studied for the first time by Taneja
\cite{tan2}, called \textit{arithmetic and geometric mean
divergence measure} or simply, \textit{AG-divergence.} More
details on these divergence measures can be seen in Taneja
\cite{tan1, tan3}.

For simplicity, we shall call the measure $D(P\vert \vert Q)$ by
\textit{relative J-divergence}, the measure $F(P\vert \vert Q)$ by
\textit{relative JS-divergence }and the measure $G(P\vert \vert
Q)$ by \textit{relative AG- divergence}. The measure $D(P\vert
\vert Q)$ is due to Dragomir et al. \cite{dgp}. The measure
$F(P\vert \vert Q)$ has been studied by many authors (Shioya and
Da-te \cite{shd}; Barnet et al. \cite{bcd} ; Lin \cite{lin}; Lin
and Wong \cite{liw}). The measure $G(P\vert \vert Q)$ we have
considered here for the first time and is a part of the measure
$T(P\vert \vert Q)$.

The one parametric generalization of the Kullback-Leibler
\cite{kul} \textit{relative information} studied in a different
way by Cressie and Read \cite{crr} is given by
\begin{equation}
\label{eq14} \Phi _s (P\vert \vert Q) = \begin{cases}
 {K_s (P\vert \vert Q) = \left[ {s(s - 1)} \right]^{ - 1}\left[
{\sum\limits_{i = 1}^n {p_i^s q_i^{1 - s} } - 1} \right],} & {s \ne 0,1} \\
 {K(Q\vert \vert P) = \sum\limits_{i = 1}^n {q_i \ln \left( {\frac{q_i }{p_i
}} \right)} ,} & {s = 0} \\
 {K(P\vert \vert Q) = \sum\limits_{i = 1}^n {p_i \ln \left( {\frac{p_i }{q_i
}} \right)} ,} & {s = 1} \\
\end{cases},
\end{equation}

\noindent for all $P,Q \in \Gamma _n $ and $s \in (-\infty, \infty)$.\\

The measure (\ref{eq14}) admits the following particular cases:
\begin{itemize}
\item[(i)] $\Phi _{ - 1} (P\vert \vert Q) = \frac{1}{2}\chi
^2(Q\vert \vert P).$

\item[(ii)] $ \Phi _0 (P\vert \vert Q) = K(Q\vert \vert P).$

\item[(iii)] $\Phi _{1 / 2} (P\vert \vert Q) = 4\left[ {1 -
B(P\vert \vert Q)} \right] = 4\,h(P\vert \vert Q)$

\item[(iv)] $\Phi _1 (P\vert \vert Q) = K(P\vert \vert Q).$

\item[(v)] $\Phi _2 (P\vert \vert Q) = \frac{1}{2}\chi ^2(P\vert
\vert Q).$
\end{itemize}

\bigskip
The measures $B(P\vert \vert Q)$, $h(P\vert \vert Q)$ and $\chi
^2(P\vert \vert Q)$ appearing in parts (iii) and (v) above, are
given by
\begin{equation}
\label{eq15}
B(P\vert \vert Q) = \sqrt {p_i q_i } ,
\end{equation}
\begin{equation}
\label{eq16}
h(P\vert \vert Q) = 1 - B(P\vert \vert Q) = \frac{1}{2}\sum\limits_{i = 1}^n
{(\sqrt {p_i } - \sqrt {q_i } )^2}
\end{equation}

\noindent and
\begin{equation}
\label{eq17}
\chi ^2(P\vert \vert Q) = \sum\limits_{i = 1}^n {\frac{(p_i - q_i )^2}{q_i
}} = \sum\limits_{i = 1}^n {\frac{p_i^2 }{q_i } - 1}
\end{equation}

\noindent respectively.

The measure $B(P\vert \vert Q)$ is famous as Bhattacharya
\cite{bha} \textit{coefficient}, the measure $h(P\vert \vert Q)$
is famous as Hellinger \cite{hel} \textit{discrimination} and the
measures $\chi ^2(P\vert \vert Q)$ is known by \textit{Chi-square}
\cite{pea} \textit{divergence}.

For more studies on the measure (\ref{eq14}) refer to Liese and
Vajda \cite{liv}, Taneja \cite{tan4}, Taneja and Kumar \cite{tak}
and Cerone et al. \cite{cdo}.

Our aim in this paper is to obtain bounds on the \textit{relative
divergence measures} that we shall classify as
\textit{non-symmetric divergence measures} given by
(\ref{eq4})-(\ref{eq7}) and on the \textit{divergence measures
}classifying as \textit{symmetric divergence measures} given by
(\ref{eq2}), (\ref{eq10}) and (\ref{eq11}) in terms of
\textit{generalized relative information} or \textit{relative
information of type s} given by (\ref{eq14}). These bounds are
studied by use of some properties of Csisz\'{a}r \cite{csi1}
\textit{f-divergence}.

\section{Csisz\'{a}r $f-$Divergence }

Given a convex function $f:[0,\infty ) \to \mathbb{R}$, the $f -
$divergence measure introduced by Csisz\'{a}r \cite{csi1} is given
by
\begin{equation}
\label{eq18}
C_f (P\vert \vert Q) =
\sum\limits_{i = 1}^n {q_i f\left( {\frac{p_i }{q_i }} \right)} ,
\end{equation}

\noindent where $P,Q \in \Gamma _n $.

It is well known in the literature \cite{csi1, csi2} that
\textit{if } $f$ \textit{ is convex and normalized, i.e., }$f(1) =
0,$ \textit{then the Csisz\'{a}r function }$C_f (P\vert \vert
Q)$\textit{ is nonnegative and convex in the pair of probability
distribution }$(P,Q) \in \Gamma _n \times \Gamma _n .$

Now we shall prove the convexity and nonnegativity of the measures
given in Section 1.

\begin{example} \label{exa21} (\textit{Relative J-divergence}). Let us
consider
\begin{equation}
\label{eq19} f_{D_1 } (x) = (x - 1)\ln \left( {\frac{x + 1}{2}}
\right), \,\, x \in (0,\infty )
\end{equation}

\noindent in (\ref{eq18}), then one gets $C_f (P\vert \vert Q) =
D(P\vert \vert Q): = D_1 $, where $D\left( {P\vert \vert Q}
\right)$ is as given by (\ref{eq4}).

Moreover,
\begin{equation}
\label{eq20}
{f}'_{D_1 } (x) = \frac{x - 1}{x + 1} + \ln \left( {\frac{x + 1}{2}}
\right)
\end{equation}

\noindent and
\begin{equation}
\label{eq21}
{f}''_{D_1 } (x) = \frac{x + 3}{(x + 1)^2}
\end{equation}

Thus from (\ref{eq21}) we see that $f_{D_1 } ^{\prime \prime }(x)
> 0$ for all $x > 0$, and hence, $f_{D_1 } (x)$ is \textit{convex}
for all $x > 0$. Also, we have $f_{D_1 } (1) = 0$. In view of this
we can say that the \textit{relative J-divergence} is
\textit{nonnegative} and \textit{convex} in the pair of
probability distributions $(P,Q) \in \Gamma _n \times \Gamma _n $.
\end{example}

\begin{example} \label{exa22} (\textit{Adjoint of relative J-divergence}).
Let us consider
\begin{equation}
\label{eq22} f_{D_2 } (x) = (1 - x)\ln \left( {\frac{x + 1}{2x}}
\right), \,\, x \in (0,\infty )
\end{equation}

\noindent in (\ref{eq18}), then one gets  $C_f (P\vert \vert Q) =
D(Q\vert \vert P): = D_2 ,$ where $D(Q\vert \vert P)$ is as given
by (\ref{eq5}).

Moreover,
\begin{equation}
\label{eq23}
{f}'_{D_2 } (x) = \frac{x - 1}{x(x + 1)} - \ln \left( {\frac{x + 1}{2x}}
\right),
\end{equation}

\noindent and
\begin{equation}
\label{eq24}
{f}''_{D_2 } (x) = \frac{3x + 1}{x^2(x + 1)^2}
\end{equation}

Thus from (\ref{eq24}) we see that ${f}''_{D_2 } (x) > 0$ for all
$x > 0$, and hence, $f_{D_2 } (x)$ is \textit{convex} for all $x >
0$. Also, we have $f_{D_2 } (1) = 0$. In view of this we can say
that the  \textit{adjoiont of relative J-divergence} is
\textit{nonnegative} and \textit{convex} in the pair of
probability distributions $(P,Q) \in \Gamma _n \times \Gamma _n $.
\end{example}

\begin{example} \label{exa23} (\textit{Relative JS-divergence}). Let us
consider
\begin{equation}
\label{eq25} f_{F_1 } (x) = \frac{1 - x}{2} - x\ln \left( {\frac{x
+ 1}{2x}} \right), \,\, x \in (0,\infty )
\end{equation}

\noindent in (\ref{eq18}), then one gets  $C_f (P\vert \vert Q) =
F(P\vert \vert Q): = F_1 $, where $F\left( {P\vert \vert Q}
\right)$ is as given by (\ref{eq6}).

Moreover,
\begin{equation}
\label{eq26}
{f}'_{F_1 } (x) = \frac{1 - x}{2(x + 1)} - \ln \left( {\frac{x + 1}{2x}}
\right),
\end{equation}

\noindent and
\begin{equation}
\label{eq27}
{f}''_{F_1 } (x) = \frac{1}{x(x + 1)^2}
\end{equation}

Thus from (\ref{eq27}) we see that ${f}''_{F_1 } (x) > 0$ for all
$x > 0$, and hence, $f_{F_1 } (x)$ is \textit{convex} for all $x >
0$. Also, we have $f_{F_1 } (1) = 0$. In view of this we can say
that the  \textit{relative JS-divergence} is \textit{nonnegative}
and \textit{convex} in the pair of probability distributions
$(P,Q) \in \Gamma _n \times \Gamma _n $.
\end{example}

\begin{example} \label{exa24} (\textit{Adjoint of relative
JS-divergence}). Let us consider
\begin{equation}
\label{eq28} f_{F_2 } (x) = \frac{x - 1}{2} - \ln \left( {\frac{x
+ 1}{2}} \right), \,\, x \in (0,\infty )
\end{equation}

\noindent in (\ref{eq18}), then one gets  $C_f (P\vert \vert Q) =
F(Q\vert \vert P): = F_2 ,$ where $F(Q\vert \vert P)$ is as given
by (\ref{eq8}).

Moreover,
\begin{equation}
\label{eq29}
{f}'_{F_2 } (x) = \frac{x - 1}{2(x + 1)}
\end{equation}

\noindent and
\begin{equation}
\label{eq30}
{f}''_{F_2 } (x) = \frac{1}{(x + 1)^2}
\end{equation}

Thus from (\ref{eq30}) we see that ${f}''_{F_2 } (x) > 0$ for all
$x > 0$, and hence, $f_{F_2 } (x)$ is \textit{convex} for all $x >
0$. Also, we have $f_{F_2 } (1) = 0$. In view of this we can say
that the  \textit{adjoint of relative JS-divergence }is
\textit{nonnegative} and \textit{convex} in the pair of
probability distributions $(P,Q) \in \Gamma _n \times \Gamma _n $.
\end{example}

\begin{example} \label{exa25} (\textit{Relative AG-divergence}). Let us
consider
\begin{equation}
\label{eq31} f_{G_1 } (x) = \frac{x - 1}{2} + \frac{x + 1}{2}\ln
\left( {\frac{x + 1}{2x}} \right), \,\, x \in (0,\infty )
\end{equation}

\noindent in (\ref{eq18}), then one gets  $C_f (P\vert \vert Q) =
G(P\vert \vert Q): = G_1 $, where $G\left( {P\vert \vert Q}
\right)$ is as given by (\ref{eq7}).

Moreover,
\begin{equation}
\label{eq32}
{f}'_{G_1 } (x) = \frac{1}{2}\left[ {\frac{x - 1}{x} + \ln \left( {\frac{x +
1}{2x}} \right)} \right],
\end{equation}

\noindent and
\begin{equation}
\label{eq33}
{f}''_{G_1 } (x) = \frac{1}{2x^2(x + 1)}
\end{equation}

Thus from (\ref{eq33}) we see that ${f}''_{G_1 } (x) > 0$ for all
$x > 0$, and hence, $f_{G_1 } (x)$ is \textit{convex} for all $x >
0$. Also, we have $f_{G_1 } (1) = 0$. In view of this we can say
that the  \textit{relative AG-divergence} is \textit{nonnegative}
and \textit{convex} in the pair of probability distributions
$(P,Q) \in \Gamma _n \times \Gamma _n $.
\end{example}

\begin{example} \label{exa26} (\textit{Adjoint of relative
AG-divergence}). Let us consider
\begin{equation}
\label{eq34} f_{G_2 } (x) = \frac{1 - x}{2} + \frac{x + 1}{2}\ln
\left( {\frac{x + 1}{2}} \right), \,\, x \in (0,\infty )
\end{equation}

\noindent in (\ref{eq18}), then one gets  $C_f (P\vert \vert Q) =
G(Q\vert \vert P): = G_2 $, where $G\left( {Q\vert \vert P}
\right)$ is as given by (\ref{eq9}).

Moreover,
\begin{equation}
\label{eq35}
{f}'_{G_2 } (x) = \frac{1}{2}\ln \left( {\frac{x + 1}{2}} \right),
\end{equation}

\noindent and
\begin{equation}
\label{eq36}
{f}''_{G_2 } (x) = \frac{1}{2(x + 1)}
\end{equation}

Thus from (\ref{eq36}) we see that ${f}''_{G_2 } (x) > 0$ for all
$x > 0$, and hence, $f_{G_2 } (x)$ is \textit{convex} for all $x >
0$. Also, we have $f_{G_2 } (1) = 0$. In view of this we can say
that the  \textit{adjoint of relative AG-divergence} is
\textit{nonnegative} and \textit{convex} in the pair of
probability distributions $(P,Q) \in \Gamma _n \times \Gamma _n $.
\end{example}

\begin{example} \label{exa27} (\textit{J-divergence}). Let us consider
\begin{equation}
\label{eq37} f_J (x) = (x - 1)\ln x, \,\, x \in (0,\infty )
\end{equation}

\noindent in (\ref{eq18}), then one gets  $C_f (P\vert \vert Q) =
J(P\vert \vert Q)$, where $J\left( {P\vert \vert Q} \right)$ is as
given by (\ref{eq2}).

Moreover,
\begin{equation}
\label{eq38}
{f}'_J (x) = 1 - x^{ - 1} + \ln x,
\end{equation}

\noindent and
\begin{equation}
\label{eq39}
{f}''_J (x) = \frac{x + 1}{x^2}
\end{equation}

Thus from (\ref{eq39}) we see that ${f}''_J (x) > 0$ for all $x >
0$, and hence, $f_J (x)$ is \textit{convex} for all $x > 0$. Also,
we have $f_J (1) = 0$. In view of this we can say that the
\textit{J-divergence }is \textit{nonnegative} and \textit{convex}
in the pair of probability distributions $(P,Q) \in \Gamma _n
\times \Gamma _n $.
\end{example}

\begin{example}\label{exa28} (\textit{JS-divergence}). Let us consider
\begin{equation}
\label{eq40} f_I (x) = \frac{x}{2}\ln x - \frac{x + 1}{2}\ln
\left( {\frac{x + 1}{2}} \right), \,\, x \in (0,\infty )
\end{equation}

\noindent in (\ref{eq18}), then one gets  $C_f (P\vert \vert Q) =
I(P\vert \vert Q)$, where $I\left( {P\vert \vert Q} \right)$ is as
given by (\ref{eq10}).

Moreover,
\begin{equation}
\label{eq41}
{f}'_I (x) = - \frac{1}{2}\ln \left( {\frac{x + 1}{2x}} \right),
\end{equation}

\noindent and
\begin{equation}
\label{eq42}
{f}''_I (x) = \frac{1}{2x(x + 1)}
\end{equation}

Thus from (\ref{eq42}) we see that ${f}''_I (x) > 0$ for all $x >
0$, and hence, $f_I (x)$ is \textit{convex} for all $x > 0$. Also,
we have $f_I (1) = 0$. In view of this we can say that the
\textit{JS-divergence }is \textit{nonnegative} and \textit{convex}
in the pair of probability distributions $(P,Q) \in \Gamma _n
\times \Gamma _n $.
\end{example}

\begin{example} \label{exa29} (\textit{AG-divergence}). Let us consider
\begin{equation}
\label{eq43} f_T (x) = \left( {\frac{x + 1}{2}} \right)\ln \left(
{\frac{x + 1}{2\sqrt x }} \right), \,\, x \in (0,\infty ),
\end{equation}

\noindent in (\ref{eq18}), then one gets $C_f (P\vert \vert Q) =
T(P\vert \vert Q),$ where $T(P\vert \vert Q)$ is as given by
(\ref{eq11}).

Moreover,
\begin{equation}
\label{eq44}
{f}'_T (x) = \frac{1}{4}\left[ {1 - x^{ - 1} + 2\ln \left( {\frac{x +
1}{2\sqrt x }} \right)} \right],
\end{equation}

\noindent and
\begin{equation}
\label{eq45}
{f}''_T (x) = \frac{1}{4}\left( {\frac{x^2 + 1}{x^3 + x^2}} \right).
\end{equation}

Thus from (\ref{eq45}) we see that ${f}''_T (x) > 0$ for all $x >
0$, and hence, $f_T (x)$ is \textit{convex} for all $x > 0$. Also,
we have $f_T (1) = 0$. In view of this we can say that the
\textit{AG-divergence} is \textit{nonnegative} and \textit{convex}
in the pair of probability distributions $(P,Q) \in \Gamma _n
\times \Gamma _n $.
\end{example}

The above examples give only the \textit{nonnegativity} and
\textit{convexity} of the \textit{symmetric} and
\textit{non-symmetric divergence measures}. Here we shall make use
of this property to get bounds in terms of \textit{relative
information of type s}. For more properties of these measures
refer to Taneja \cite{tan7}.

\section{Csisz\'{a}r $f-$Divergence and Relative
Information of Type $s$}

The following two theorems are due to Taneja \cite{tan4} and
Taneja and Kumar \cite{tak}.

\begin{theorem} \label{the31} Let $P,Q \in \Gamma _n $ and $s \in \mathbb{R}: = ( -
\infty ,\infty )$, then we have
\begin{equation}
\label{eq46} 0 \leqslant \Phi _s (P\vert \vert Q) \leqslant
E_{\Phi _s } (P\vert \vert Q),
\end{equation}

\noindent where
\begin{equation}
\label{eq47} E_{\Phi _s} (P\vert \vert Q) = \begin{cases}
 {(s - 1)^{ - 1}\sum\limits_{i = 1}^n {(p_i - q_i )\left( {\frac{p_i }{q_i
}} \right)^{s - 1},} } & {s \ne 1} \\
 {\sum\limits_{i = 1}^n {(p_i - q_i )\ln \left( {\frac{p_i }{q_i }}
\right),} } & {s = 1} \\
\end{cases}.
\end{equation}

Let $P,Q \in \Gamma _n $ be such that there exists $r,R$ with $0 <
r \leqslant \frac{p_i }{q_i } \leqslant R < \infty $, $\forall i
\in \{1,2,...,n\}$, then
\begin{equation}
\label{eq48} 0 \leqslant \Phi _s (P\vert \vert Q) \leqslant
A_{\Phi _s } (r,R)
\end{equation}

\noindent where
\begin{equation}
\label{eq49} A_{\Phi _s } (r,R) = \frac{1}{4}(R - r)^2\mbox{
}\begin{cases}
 {\frac{R^{s - 1} - r^{s - 1}}{(R - r)(s - 1)},} & {s \ne 1} \\
 {\frac{\ln R - \ln r}{R - r},} & {s = 1} \\
\end{cases}
\end{equation}

Futher, if we suppose that $0 < r \leqslant 1 \leqslant R < \infty
$, $r \ne R$, then
\begin{equation}
\label{eq50} 0 \leqslant \Phi _s (P\vert \vert Q) \leqslant
B_{\Phi _s} (r,R)
\end{equation}

\noindent where
\begin{equation}
\label{eq51} B_{\Phi _s } (r,R) = \begin{cases}
 {\frac{(R - 1)(r^s - 1) + (1 - r)(R^s - 1)}{(R - r)s(s - 1)},} & {s \ne
0,1} \\
 {\frac{(R - 1)\ln \frac{1}{r} + (1 - r)\ln \frac{1}{R}}{(R - r)},} & {s =
0} \\
 {\frac{(R - 1)r\ln r + (1 - r)R\ln R}{(R - r)},} & {s = 1} \\
\end{cases}
\end{equation}

Moreover, the following inequalities hold:
\begin{equation}
\label{eq52} E_{\Phi _s } (P\vert \vert Q) \leqslant A_{\Phi _s }
(r,R),
\end{equation}
\begin{equation}
\label{eq53} B_{\Phi _s } (r,R) \leqslant A_{\Phi _s } (r,R)
\end{equation}

\noindent and
\begin{equation}
\label{eq54} B_{\Phi _s } (r,R) - \Phi _s (P\vert \vert Q)
\leqslant A_{\Phi _s } (r,R).
\end{equation}
\end{theorem}

\begin{theorem} \label{the32}  Let $f:I \subset \mathbb{R}_ + \to
\mathbb{R}$ the generating mapping is normalized, i.e., $f(1) = 0$
and satisfy the assumptions:
\begin{itemize}

\item[(i)] $f$is twice differentiable on $(r,R)$;

\item[(ii)] there exists real constants $m,M$ such that $0 < m <
M$ and
\end{itemize}
\begin{equation}
\label{eq55} m \leqslant x^{2 - s}{f}''(x) \leqslant M, \quad
\forall x \in (r,R), \quad
 - \infty < s < \infty .
\end{equation}

\noindent then, we have
\begin{equation}
\label{eq56} m \Phi _s (P\vert \vert Q)
 \leqslant C_f (P\vert \vert Q) \leqslant M
\Phi _s (P\vert \vert Q)
\end{equation}

\noindent and
\begin{align}
\label{eq57} m\left[ {E_{\Phi _s } (P\vert \vert Q) - \Phi _s
(P\vert \vert Q)} \right] & \leqslant E_{C_f } (P\vert \vert Q) -
C_f (P\vert \vert Q)\\
& \leqslant M\left[ {E_{\Phi _s } (P\vert \vert Q) - \Phi _s
(P\vert \vert Q)} \right].\notag
\end{align}

Let $P,Q \in \Gamma _n $ be such that there exists $r,R$ with $0 <
r \leqslant \frac{p_i }{q_i } \leqslant R < \infty $, $\forall i
\in \{1,2,...,n\}$, then
\begin{align}
\label{eq58} m\left[ {A_{\Phi _s } (r,R) - \Phi _s (P\vert \vert
Q)} \right] & \leqslant A_{C_f } (r,R) - C_f (P\vert \vert Q)\\
& \leqslant M\left[ {A_{\Phi _s } (r,R) - \Phi _s (P\vert \vert
Q)} \right]\notag
\end{align}

Further, if we suppose that $0 < r \leqslant 1 \leqslant R <
\infty $, $r \ne R$, then
\begin{align}
\label{eq59} m\left[ {B_{\Phi _s } (r,R) - \Phi _s (P\vert \vert
Q)} \right] & \leqslant B_{C_f } (r,R) - C_f (P\vert \vert Q)\\
& \leqslant M\left[ {B_{\Phi _s } (r,R) - \Phi _s (P\vert \vert
Q)} \right].\notag
\end{align}
\end{theorem}

The Theorem \ref{the31} is obtained by applying some of the
results due to Dragomir \cite{dra1, dra2}. The Theorem \ref{the32}
unifies some of the results studied by Dragomir \cite{dra3, dra4,
dra5}. For an improved version of Theorem \ref{the32} refers to
Taneja \cite{tan5}.

The aim here is to apply Theorem \ref{the32} by taking different
values of $f$ given by examples \ref{exa21}-\ref{exa29}. This we
have done only applying the inequalities (\ref{eq56}), while the
results for the inequalities (\ref{eq57})-(\ref{eq59}) can be done
on similar lines. These details are omitted here.

\section{Bounds On Non-Symmetric Divergence Measures}

In this section, we have applied the inequalities (\ref{eq56}) and
used the condition (\ref{eq55}) to obtain bounds for the measures
given in (\ref{eq4})-(\ref{eq9}).

\begin{theorem} \label{the41} The following bounds on
\textit{relative J-divergence} hold:
\begin{equation}
\label{eq60}
\frac{r^{2 - s}(r + 3)}{(r + 1)^2}\Phi _s (P\vert \vert Q)
 \leqslant D(P\vert \vert Q) \leqslant
\frac{R^{2 - s}(R + 3)}{(R + 1)^2}\Phi _s (P\vert \vert Q), \,\, s
\leqslant \frac{3}{4}
\end{equation}

\noindent and
\begin{equation}
\label{eq61}
\frac{R^{2 - s}(R + 3)}{(R + 1)^2}\Phi _s (P\vert \vert Q)
 \leqslant D(P\vert \vert Q) \leqslant
\frac{r^{2 - s}(r + 3)}{(r + 1)^2}\Phi _s (P\vert \vert Q), \,\, s
\geqslant 2.
\end{equation}
\end{theorem}

\begin{proof} Let us consider
\begin{equation}
\label{eq62} g_{D_1 } (x) = x^{2 - s}{f}''_{D_1 } (x) = \frac{x^{2
- s}(x + 3)}{(x + 1)^2}, \,\, x \in (0,\infty ),
\end{equation}

\noindent where ${f}''_{D_1 } (x)$ is as given by (\ref{eq21}).

From (\ref{eq62}), one can get
\begin{equation}
\label{eq63}
{g}'_{D_1 } (x) = - \frac{x^{1 - s}\left[ {(s - 1)x^2 + (4s - 3)x + 3(s -
2)} \right]}{(x + 1)^3}
\begin{cases}
 { \geqslant 0,} & {s \leqslant \tfrac{3}{4}} \\
 { \leqslant 0,} & {s \geqslant 2} \\
\end{cases}.
\end{equation}

In view of (\ref{eq63}), we conclude that
\begin{equation}
\label{eq64} m = \mathop {\inf }\limits_{x \in [r,R]} g_{D_1 } (x)
= \begin{cases}
 {\frac{r^{2 - s}(r + 3)}{(r + 1)^2},} & {s \leqslant \tfrac{3}{4}} \\
 {\frac{R^{2 - s}(R + 3)}{(R + 1)^2},} & {s \geqslant 2} \\
\end{cases}
\end{equation}

\noindent and
\begin{equation}
\label{eq65} M = \mathop {\sup }\limits_{x \in [r,R]} g_{D_1 } (x)
= \begin{cases}
 {\frac{R^{2 - s}(R + 3)}{(R + 1)^2},} & {s \leqslant \tfrac{3}{4}} \\
 {\frac{r^{2 - s}(r + 3)}{(r + 1)^2},} & {s \geqslant 2} \\
\end{cases}.
\end{equation}

In view of (\ref{eq64}), (\ref{eq65}) and (\ref{eq56}), we get the
inequalities (\ref{eq60}) and (\ref{eq61}).
\end{proof}

Some particular cases of the Theorem \ref{the41} are summarized in
the following corollary.

\begin{corollary} \label{cor41} The following bounds hold:
\begin{equation}
\label{eq66}
\frac{r^3(r + 3)}{2(r + 1)^2}\chi ^2(Q\vert \vert P) \leqslant D(P\vert
\vert Q) \leqslant \frac{R^3(R + 3)}{2(R + 1)^2}\chi ^2(Q\vert \vert P),
\end{equation}
\begin{equation}
\label{eq67}
\frac{r^2(r + 3)}{(r + 1)^2}K(Q\vert \vert P) \leqslant D(P\vert \vert Q)
\leqslant \frac{R^2(R + 3)}{(R + 1)^2}K(Q\vert \vert P),
\end{equation}
\begin{equation}
\label{eq68}
\frac{4r^{3 / 2}(r + 3)}{(r + 1)^2}h(P\vert \vert Q) \leqslant D(P\vert
\vert Q) \leqslant \frac{4R^{3 / 2}(R + 3)}{(R + 1)^2}h(P\vert \vert Q)
\end{equation}

\noindent and
\begin{equation}
\label{eq69}
\frac{R + 3}{2(R + 1)^2}\chi ^2(P\vert \vert Q) \leqslant D(P\vert \vert Q)
\leqslant \frac{r + 3}{2(r + 1)^2}\chi ^2(P\vert \vert Q).
\end{equation}
\end{corollary}

\begin{proof} Inequalities (\ref{eq66}), (\ref{eq67}) and
(\ref{eq68}) follows from (\ref{eq60}) by taking $s = - 1$, $s =
0$ and $s = \frac{1}{2}$ respectively. The inequalities
(\ref{eq69}) follow form (\ref{eq61}) by taking $s = 2$.
\end{proof}

The case $s = 1$ is not included in the inequalities (\ref{eq60})
and (\ref{eq61}). This we shall do separately in the following
proposition.

\begin{proposition} \label{pro41} The following inequality hold:
\begin{equation} \label{eq70} D(P\vert \vert Q) \leqslant
\frac{9}{8}K(P\vert \vert Q).
\end{equation}
\end{proposition}

\begin{proof} For $s = 1$ in (\ref{eq62}), we have
\begin{equation}
\label{eq71}
g_D (x) = \frac{x(x + 3)}{(x + 1)^2}.
\end{equation}

This gives
\begin{equation}
\label{eq72}
{g}'_D (x) = - \frac{x - 3}{(x + 1)^3}
\begin{cases}
 { \geqslant 0,} & {x \leqslant 3} \\
 { \leqslant 0,} & {x \geqslant 3} \\
\end{cases}.
\end{equation}

Thus we conclude from (\ref{eq72}) that the function $g_D (x)$
given by (\ref{eq71}) is increasing in $x \in (0,3)$ and
decreasing in $x \in (3,\infty )$, and hence
\begin{equation}
\label{eq73} M = \mathop {\sup }\limits_{x \in (0,\infty )} g_D
(x) = \frac{9}{8}.
\end{equation}

Now (\ref{eq73}) together with (\ref{eq56}) give the required
result.
\end{proof}

\begin{theorem} \label{the42} The following bounds on \textit{adjoint
}of\textit{ relative J-divergence} holds:
\begin{equation}
\label{eq74}
\frac{r^{ - s}(3r + 1)}{(r + 1)^2}\Phi _s (P\vert \vert Q)
 \leqslant D(Q\vert \vert P) \leqslant
\frac{R^{ - s}(3R + 1)}{(R + 1)^2}\Phi _s (P\vert \vert Q), \,\, s
\leqslant - 1
\end{equation}

\noindent and
\begin{equation}
\label{eq75}
\frac{R^{ - s}(3R + 1)}{(R + 1)^2}\Phi _s (P\vert \vert Q)
 \leqslant D(Q\vert \vert P) \leqslant
\frac{r^{ - s}(3r + 1)}{(r + 1)^2}\Phi _s (P\vert \vert Q), \,\, s
\geqslant \frac{1}{4}
\end{equation}
\end{theorem}

\begin{proof} Let us consider
\begin{equation}
\label{eq76} g_{D_2 } (x) = x^{2 - s}{f}''_{D_2 } (x) = \frac{x^{
- s}(3x + 1)}{(x + 1)^2}, \,\, x \in (0,\infty ),
\end{equation}

\noindent where ${f}''_{D_2 } (x)$ is as given by (\ref{eq24}).

From (\ref{eq76}), one can get
\begin{equation}
\label{eq77}
{g}'_D (x) = - \frac{x^{ - s}\left[ {3(s + 1)x^2 + (4s - 1)x + s}
\right]}{(x + 1)^3}
\begin{cases}
 { \geqslant 0,} & {s \leqslant - 1} \\
 { \leqslant 0,} & {s \geqslant \frac{1}{4}} \\
\end{cases}.
\end{equation}

In view of (\ref{eq77}), we conclude that
\begin{equation}
\label{eq78} m = \mathop {\inf }\limits_{x \in [r,R]} g_D (x) =
\mathop {\min }\limits_{x \in [r,R]} g_D (x) = \begin{cases}
 {\frac{r^{ - s}(3r + 1)}{(r + 1)^2},} & {s \leqslant - 1} \\
 {\frac{R^{ - s}(3R + 1)}{(R + 1)^2},} & {s \geqslant \frac{1}{4}} \\
\end{cases}
\end{equation}

\noindent and
\begin{equation}
\label{eq79} M = \mathop {\sup }\limits_{x \in [r,R]} g_D (x) =
\begin{cases}
 {\frac{R^{ - s}(3R + 1)}{(R + 1)^2},} & {s \leqslant - 1} \\
 {\frac{r^{ - s}(3r + 1)}{(r + 1)^2},} & {s \geqslant \frac{1}{4}} \\
\end{cases}.
\end{equation}

Now from (\ref{eq78}), (\ref{eq79}) and (\ref{eq56}), we get the
inequalities (\ref{eq74}) and (\ref{eq75}).
\end{proof}

Some particular cases of the Theorem \ref{the42} are summarized in
the following corollary.

\begin{corollary} \label{cor42} The following bounds hold:
\begin{equation}
\label{eq80}
\frac{r(3r + 1)}{(r + 1)^2}\chi ^2(Q\vert \vert P) \leqslant D(Q\vert \vert
P) \leqslant \frac{R(3R + 1)}{(R + 1)^2}\chi ^2(Q\vert \vert P),
\end{equation}
\begin{equation}
\label{eq81}
\frac{4(3R + 1)}{\sqrt R (R + 1)^2}h(P\vert \vert Q) \leqslant D(Q\vert
\vert P) \leqslant \frac{4(3r + 1)}{\sqrt r (r + 1)^2}h(P\vert \vert Q),
\end{equation}
\begin{equation}
\label{eq82}
\frac{3R + 1}{R(R + 1)^2}K(P\vert \vert Q) \leqslant D(Q\vert \vert P)
\leqslant \frac{3R + 1}{R(R + 1)^2}K(P\vert \vert Q)
\end{equation}

\noindent and
\begin{equation}
\label{eq83}
\frac{3R + 1}{R^2(R + 1)^2}\chi ^2(P\vert \vert Q) \leqslant D(Q\vert \vert
P) \leqslant \frac{3r + 1}{r^2(r + 1)^2}\chi ^2(P\vert \vert Q).
\end{equation}
\end{corollary}

\begin{proof} Inequalities (\ref{eq80}) follows from
(\ref{eq74}) by taking $s = - 1$. The inequalities (\ref{eq81}),
(\ref{eq82}) and (\ref{eq83}) follows from (\ref{eq75}) by taking
$s = \frac{1}{2}$, $s = 1$ and $s = 2$ respectively.
\end{proof}

The case $s = 0$ is not included in the inequalities (\ref{eq74})
and (\ref{eq75}). This we shall do separately in the following
proposition.

\begin{proposition} \label{pro42} The following bound hold:
\begin{equation}
\label{eq84} D(Q\vert \vert P) \leqslant \frac{9}{8}K(Q\vert \vert
P),
\end{equation}
\end{proposition}

\begin{proof} For $s = 0$ in (\ref{eq76}), we have
\begin{equation}
\label{eq85}
g_{D_2 } (x) = \frac{3x + 1}{(x + 1)^2}.
\end{equation}

This gives
\begin{equation}
\label{eq86}
{g}'_{D_2 } (x) = - \frac{3x - 1}{(x + 1)^3}
\begin{cases}
 { \geqslant 0,} & {x \leqslant \frac{1}{3}} \\
 { \leqslant 0,} & {x \geqslant \frac{1}{3}} \\
\end{cases}.
\end{equation}

Thus from (\ref{eq86}) we conclude that the function $g_{D_2 }
(x)$ given by (\ref{eq85}) is increasing in $x \in
(0,\frac{1}{3})$ and decreasing in $x \in (\frac{1}{3},\infty )$,
and hence
\begin{equation}
\label{eq87}
M = \mathop {\sup }\limits_{x \in (0,\infty )} g_{D_2 } (x) = g_{D_2 }
(\frac{1}{3}) = \frac{9}{8}.
\end{equation}

Now (\ref{eq87}) together with (\ref{eq56}) give the required
result.
\end{proof}

\begin{theorem} \label{the43} The following bounds on \textit{relative
JS-divergence} hold:
\begin{equation}
\label{eq88}
\frac{r^{1 - s}}{(r + 1)^2}\Phi _s (P\vert \vert Q)
 \leqslant F(P\vert \vert Q) \leqslant
\frac{R^{1 - s}}{(R + 1)^2}\Phi _s (P\vert \vert Q), \,\, s
\leqslant - 1
\end{equation}

\noindent and
\begin{equation}
\label{eq89}
\frac{R^{1 - s}}{(R + 1)^2}\Phi _s (P\vert \vert Q)
 \leqslant F(P\vert \vert Q) \leqslant
\frac{r^{1 - s}}{(r + 1)^2}\Phi _s (P\vert \vert Q), \,\, s
\geqslant 1.
\end{equation}
\end{theorem}

\begin{proof} Let us consider
\begin{equation}
\label{eq90} g_{F_1 } (x) = x^{2 - s}{f}''_{F_1 } (x) = \frac{x^{1
- s}}{(x + 1)^2}, \,\, x \in (0,\infty )
\end{equation}

\noindent where ${f}''_{F_1 } (x)$ is as given by (\ref{eq27}).

From (\ref{eq90}) one can get
\begin{equation}
\label{eq91}
{g}'_{F_1 } (x) = - \frac{x^{ - s}\left[ {(s + 1)x + (s - 1)} \right]}{(x +
1)^3}
\begin{cases}
 { \geqslant 0,} & {s \leqslant - 1} \\
 { \leqslant 0,} & {s \geqslant 1} \\
\end{cases}.
\end{equation}

In view of (\ref{eq91}), we conclude that
\begin{equation}
\label{eq92} m = \mathop {\inf }\limits_{x \in [r,R]} g_{F_1 } (x)
= \begin{cases}
 {\frac{r^{1 - s}}{(r + 1)^2},} & {s \leqslant - 1} \\
 {\frac{R^{1 - s}}{(R + 1)^2},} & {s \geqslant 1} \\
\end{cases}
\end{equation}

\noindent and
\begin{equation}
\label{eq93} M = \mathop {\sup }\limits_{x \in [r,R]} g_{F_1 } (x)
= \begin{cases}
 {\frac{R^{1 - s}}{(R + 1)^2},} & {s \leqslant - 1} \\
 {\frac{r^{1 - s}}{(r + 1)^2},} & {s \geqslant 1} \\
\end{cases}.
\end{equation}

Now (\ref{eq92}) and (\ref{eq93}) together with (\ref{eq56}) give
the required result.
\end{proof}

Some particular cases of the Theorem \ref{the43} are summarized in
the following corollary.

\begin{corollary} \label{cor43} The following bounds hold:
\begin{equation}
\label{eq94}
\frac{r^2}{2(r + 1)^2}\chi ^2(Q\vert \vert P) \leqslant F(P\vert \vert Q)
\leqslant \frac{R^2}{2(R + 1)^2}\chi ^2(Q\vert \vert P),
\end{equation}
\begin{equation}
\label{eq95}
\frac{1}{(R + 1)^2}K(P\vert \vert Q) \leqslant F(P\vert \vert Q) \leqslant
\frac{1}{(r + 1)^2}K(P\vert \vert Q)
\end{equation}

\noindent and
\begin{equation}
\label{eq96}
\frac{1}{2R(R + 1)^2}\chi ^2(P\vert \vert Q) \leqslant F(P\vert \vert Q)
\leqslant \frac{1}{2r(r + 1)^2}\chi ^2(P\vert \vert Q).
\end{equation}
\end{corollary}

\begin{proof} Inequalities (\ref{eq94}) follows from
(\ref{eq88}) by taking $s = - 1$. The inequalities (\ref{eq95})
and (\ref{eq96}) follows from (\ref{eq89}) by taking $s = 1$ and
$s = 2$ respectively.
\end{proof}

The cases $s = 0$ and $s = \frac{1}{2}$ are not included in the
inequalities (\ref{eq88}) and (\ref{eq89}). This we shall do
separately in the following proposition.

\begin{proposition} \label{pro43} The following bounds hold:
\begin{equation}
\label{eq97}
F(P\vert \vert Q) \leqslant \frac{1}{4}K(Q\vert \vert P)
\end{equation}

\noindent and
\begin{equation}
\label{eq98}
F(P\vert \vert Q) \leqslant \frac{3\sqrt 3 }{4}h(P\vert \vert Q).
\end{equation}
\end{proposition}

\begin{proof} For $s = 0$ in (\ref{eq90}), we have
\begin{equation}
\label{eq99}
g_{F_1 } (x) = \frac{x}{(x + 1)^2}.
\end{equation}

This gives
\begin{equation}
\label{eq100}
{g}'_{F_1 } (x) = - \frac{x - 1}{(x + 1)^3}
\begin{cases}
 { \geqslant 0,} & {x \leqslant 1} \\
 { \leqslant 0,} & {x \geqslant 1} \\
\end{cases}.
\end{equation}

From (\ref{eq100}) we conclude that the function $g_{F_1 } (x)$
given by (\ref{eq100}) is increasing in $x \in (0,1)$ and
decreasing in $x \in (1,\infty )$, and hence
\begin{equation}
\label{eq101}
M = \mathop {\sup }\limits_{x \in (0,\infty )} g_{F_1 } (x) = g_{F_1 } (1) =
\frac{1}{4}.
\end{equation}

Now (\ref{eq101}) together with (\ref{eq56}) give the inequality
(\ref{eq97}).

Again let us take $s = \frac{1}{2}$ in (\ref{eq90}), we have
\begin{equation}
\label{eq102}
g_{F_1 } (x) = \frac{\sqrt x }{(x + 1)^2}.
\end{equation}

This gives
\begin{equation}
\label{eq103}
{g}'_{F_1 } (x) = - \frac{3x - 1}{2\sqrt x (x + 1)^3}
\begin{cases}
 { \geqslant 0,} & {x \leqslant \frac{1}{3}} \\
 { \leqslant 0,} & {x \geqslant \frac{1}{3}} \\
\end{cases}.
\end{equation}

Thus from (\ref{eq103}), we conclude that the function $g_{F_1 }
(x)$ given by (\ref{eq102}) is increasing in $x \in
(0,\frac{1}{3})$ and decreasing in $x \in (\frac{1}{3},\infty )$,
and hence
\begin{equation}
\label{eq104}
M = \mathop {\sup }\limits_{x \in (0,\infty )} g_{F_1 } (x) = g_{F_1 }
(\frac{1}{3}) = \frac{3\sqrt 3 }{16}.
\end{equation}

Now (\ref{eq104}) together with (\ref{eq56}) give the inequalities
(\ref{eq98}).
\end{proof}

\begin{theorem} \label{the44} The following bounds on \textit{adjoint
} of \textit{relative JS-divergence} hold:
\begin{equation}
\label{eq105}
\frac{r^{2 - s}}{(r + 1)^2}\Phi _s (P\vert \vert Q)
 \leqslant F(Q\vert \vert P) \leqslant
\frac{R^{2 - s}}{(R + 1)^2}\Phi _s (P\vert \vert Q), \,\, s
\leqslant 0
\end{equation}

\noindent and
\begin{equation}
\label{eq106}
\frac{R^{2 - s}}{(R + 1)^2}\Phi _s (P\vert \vert Q)
 \leqslant F(Q\vert \vert P) \leqslant
\frac{r^{2 - s}}{(r + 1)^2}\Phi _s (P\vert \vert Q), \,\, s
\geqslant 2
\end{equation}
\end{theorem}

\begin{proof} Let us consider
\begin{equation}
\label{eq107} g_{F_2 } (x) = x^{2 - s}{f}''_{F_2 } (x) =
\frac{x^{2 - s}}{(x + 1)^2}, \,\, x \in (0,\infty )
\end{equation}

\noindent where ${f}''_{F_2 } (x)$ is as given by (\ref{eq30}).

From (\ref{eq107}) one can get
\begin{equation}
\label{eq108} {g}'_{F_2 } (x) = - \frac{x^{1 - s}\left[ {sx + (s -
2)} \right]}{(x + 1)^3}
\begin{cases}
 { \geqslant 0,} & {s \leqslant 0} \\
 { \leqslant 0,} & {s \geqslant 2} \\
\end{cases}.
\end{equation}

In view of (\ref{eq108}), we conclude that
\begin{equation}
\label{eq109} m = \mathop {\inf }\limits_{x \in [r,R]} g_{F_2} (x)
=
\begin{cases}
 {\frac{r^{2 - s}}{(r + 1)^2},} & {s \leqslant 0} \\
 {\frac{R^{2 - s}}{(R + 1)^2},} & {s \geqslant 2} \\
\end{cases}
\end{equation}

\noindent and
\begin{equation}
\label{eq110} M = \mathop {\sup }\limits_{x \in [r,R]} g_{F_2} (x)
=
\begin{cases}
 {\frac{R^{2 - s}}{(R + 1)^2},} & {s \leqslant 0} \\
 {\frac{r^{2 - s}}{(r + 1)^2},} & {s \geqslant 2} \\
\end{cases}.
\end{equation}

Now (\ref{eq109}) and (\ref{eq110}) together with (\ref{eq56})
give the required result.
\end{proof}

Some particular cases of the Theorem \ref{the44} are summarized in
the following corollary.

\begin{corollary} \label{cor44} The following bounds hold:
\begin{equation}
\label{eq111}
\frac{r^3}{2(r + 1)^2}\chi ^2(Q\vert \vert P) \leqslant F(Q\vert \vert P)
\leqslant \frac{R^3}{2(R + 1)^2}\chi ^2(Q\vert \vert P),
\end{equation}
\begin{equation}
\label{eq112}
\frac{r^2}{(r + 1)^2}K(Q\vert \vert P) \leqslant F(Q\vert \vert P) \leqslant
\frac{R^2}{(R + 1)^2}K(Q\vert \vert P)
\end{equation}

\noindent and
\begin{equation}
\label{eq113}
\frac{1}{(R + 1)^2}\chi ^2(P\vert \vert Q) \leqslant F(Q\vert \vert P)
\leqslant \frac{1}{(r + 1)^2}\chi ^2(P\vert \vert Q),
\end{equation}
\end{corollary}

\begin{proof} Inequalities (\ref{eq111}) and (\ref{eq112})
follows from (\ref{eq105}) by taking $s = - 1$ and $s = 0$
respectively. The inequalities (\ref{eq113}) follow from
(\ref{eq106}) by taking $s = 2$.
\end{proof}

The cases $s = \frac{1}{2}$ and $s = 1$ are not included in the
inequalities (\ref{eq105}) and (\ref{eq106}). This we shall do
separately in the following proposition.

\begin{proposition} \label{pro44} The following bounds hold:
\begin{equation}
\label{eq114}
F(Q\vert \vert P) \leqslant \frac{3\sqrt 3 }{4}h(P\vert \vert Q),
\end{equation}

\noindent and
\begin{equation}
\label{eq115} F(Q\vert \vert P) \leqslant \frac{1}{4}K(P\vert
\vert Q).
\end{equation}
\end{proposition}

\begin{proof} For $s = \frac{1}{2}$ in (\ref{eq107}), we
have
\begin{equation}
\label{eq116}
g_{F_2 } (x) = \frac{x\sqrt x }{(x + 1)^2}.
\end{equation}

This gives
\begin{equation}
\label{eq117}
{g}'_{F_2 } (x) = - \frac{\sqrt x (x - 3)}{2(x + 1)^3}
\begin{cases}
 { \geqslant 0,} & {x \leqslant 3} \\
 { \leqslant 0,} & {x \geqslant 3} \\
\end{cases}.
\end{equation}

Thus from (\ref{eq117}) we conclude that the function $g_{F_2 }
(x)$ given by (\ref{eq116}) is increasing in $x \in (0,3)$ and
decreasing in $x \in (3,\infty )$, and hence
\begin{equation}
\label{eq118}
M = \mathop {\sup }\limits_{x \in (0,\infty )} g_{F_2 } (x) = g_{F_2 } (3) =
\frac{3\sqrt 3 }{16}.
\end{equation}

Now (\ref{eq118}) together with (\ref{eq56}) give the inequality
(\ref{eq114}).

Again for $s = 0$ in (\ref{eq107}), we have
\begin{equation}
\label{eq119}
g_{F_2 } (x) = \frac{x}{(x + 1)^2}.
\end{equation}

This gives
\begin{equation}
\label{eq120}
{g}'_{F_2 } (x) = - \frac{x - 1}{(x + 1)^3}
\begin{cases}
 { \geqslant 0,} & {x \leqslant 1} \\
 { \leqslant 0,} & {x \geqslant 1} \\
\end{cases}.
\end{equation}

Thus from (\ref{eq120}), we conclude that the function $g_{F_2 }
(x)$ given by (\ref{eq119}) is increasing in $x \in (0,1)$ and
decreasing in $x \in (1,\infty )$, and hence
\begin{equation}
\label{eq121}
M = \mathop {\sup }\limits_{x \in (0,\infty )} g_{F_2 } (x) = g_{F_2 } (1) =
\frac{1}{4}.
\end{equation}

Now (\ref{eq121}) together with (\ref{eq56}) give the inequality
(\ref{eq115}).
\end{proof}

\begin{theorem} \label{the45} The following bounds on \textit{relative
AG-divergence} hold:
\begin{equation}
\label{eq122}
\frac{1}{2r^s(r + 1)}\Phi _s (P\vert \vert Q)
 \leqslant G(P\vert \vert Q) \leqslant
\frac{1}{2R^s(R + 1)}\Phi _s (P\vert \vert Q), \,\, s \leqslant -
1
\end{equation}

\noindent and
\begin{equation}
\label{eq123}
\frac{1}{2R^s(R + 1)}\Phi _s (P\vert \vert Q)
 \leqslant G(P\vert \vert Q) \leqslant
\frac{r^{2 - s}(r + 3)}{(r + 1)^2}\Phi _s (P\vert \vert Q), \,\, s
\geqslant 0
\end{equation}
\end{theorem}

\begin{proof} Let us consider
\begin{equation}
\label{eq124} g_{G_1 } (x) = x^{2 - s}{f}''_{G_1 } (x) = \frac{x^{
- s}}{2(x + 1)}, \,\, x \in (0,\infty ),
\end{equation}

\noindent where ${f}''_{G_1 } (x)$ is as given by (\ref{eq33}).

From (\ref{eq124}) one gets
\begin{equation}
\label{eq125}
{g}'_{G_1 } (x) = - \frac{x^{ - 1 - s}\left[ {(s + 1)x + s} \right]}{2(x +
1)^2}
\begin{cases}
 { \geqslant 0,} & {s \leqslant - 1} \\
 { \leqslant 0,} & {s \geqslant 0} \\
\end{cases}.
\end{equation}

In view of (\ref{eq125}), we conclude that
\begin{equation}
\label{eq126} m = \mathop {\inf }\limits_{x \in [r,R]} g_{G_1} (x)
= \begin{cases}
 {\frac{1}{2r^s(r + 1)},} & {s \leqslant - 1} \\
 {\frac{1}{2R^s(R + 1)},} & {s \geqslant 0} \\
\end{cases}
\end{equation}

\noindent and
\begin{equation}
\label{eq127} M = \mathop {\sup }\limits_{x \in [r,R]} g_{G_1} (x)
=
\begin{cases}
 {\frac{1}{2R^s(R + 1)},} & {s \leqslant - 1} \\
 {\frac{1}{2r^s(r + 1)},} & {s \geqslant 0} \\
\end{cases}.
\end{equation}

Now (\ref{eq126}) and (\ref{eq127}) together with (\ref{eq56})
give the required result.
\end{proof}

Some particular cases of the Theorem \ref{the45} are summarized in
the following corollary.

\begin{corollary} \label{cor45} The following bounds hold:
\begin{equation}
\label{eq128}
\frac{r}{4(r + 1)}\chi ^2(Q\vert \vert P) \leqslant G(P\vert \vert Q)
\leqslant \frac{R}{4(R + 1)}\chi ^2(Q\vert \vert P),
\end{equation}
\begin{equation}
\label{eq129}
\frac{1}{2(R + 1)}K(Q\vert \vert P) \leqslant G(P\vert \vert Q) \leqslant
\frac{1}{2(r + 1)}K(Q\vert \vert P),
\end{equation}
\begin{equation}
\label{eq130}
\frac{2}{\sqrt R (R + 1)}h(P\vert \vert Q) \leqslant G(P\vert \vert Q)
\leqslant \frac{2}{\sqrt r (r + 1)}h(P\vert \vert Q),
\end{equation}
\begin{equation}
\label{eq131}
\frac{1}{2R(R + 1)}K(P\vert \vert Q) \leqslant G(P\vert \vert Q) \leqslant
\frac{1}{2r(r + 1)}K(P\vert \vert Q),
\end{equation}

\noindent and
\begin{equation}
\label{eq132}
\frac{1}{4R^2(R + 1)}\chi ^2(P\vert \vert Q) \leqslant G(P\vert \vert Q)
\leqslant \frac{1}{4r^2(r + 1)}\chi ^2(P\vert \vert Q).
\end{equation}
\end{corollary}

\begin{proof} Inequalities (\ref{eq128}) follows from
(\ref{eq122}) by taking $s = - 1$. The inequalities (\ref{eq129}),
(\ref{eq130}), (\ref{eq131}) and (\ref{eq132}) follows from
(\ref{eq123}) by taking $s = 0$, $s = \frac{1}{2}$, $s = 1$ and $s
= 2$ respectively.
\end{proof}

\begin{theorem} \label{the46} The following bounds on \textit{adjoint
}of\textit{ relative AG-divergence} holds:
\begin{equation}
\label{eq133}
\frac{r^2}{2r^s(r + 1)}\Phi _s (P\vert \vert Q)
 \leqslant G(Q\vert \vert P) \leqslant
\frac{R^2}{2R^s(R + 1)}\Phi _s (P\vert \vert Q), \,\, s \leqslant
1
\end{equation}

\noindent and
\begin{equation}
\label{eq134}
\frac{R^2}{2R^s(R + 1)}\Phi _s (P\vert \vert Q)
 \leqslant G(Q\vert \vert P) \leqslant
\frac{r^2}{2r^s(r + 1)}\Phi _s (P\vert \vert Q), \,\, s \geqslant
2
\end{equation}
\end{theorem}

\begin{proof} Let us consider
\begin{equation}
\label{eq135} g_{G_2 } (x) = x^{2 - s}{f}''_{G_2 } (x) =
\frac{x^{2 - s}}{2(x + 1)}, \,\, x \in (0,\infty ),
\end{equation}

\noindent where ${f}''_{G_2 } (x)$ is as given by (\ref{eq36}).

From (\ref{eq135}) one gets
\begin{equation}
\label{eq136}
{g}'_{G_2 } (x) = - \frac{x^{1 - s}\left[ {(s - 1)x + (s - 2)} \right]}{2(x
+ 1)^2}
\begin{cases}
 { \geqslant 0,} & {s \leqslant 1} \\
 { \leqslant 0,} & {s \geqslant 2} \\
\end{cases}.
\end{equation}

In view of (\ref{eq136}) we conclude that
\begin{equation}
\label{eq137} m = \mathop {\inf }\limits_{x \in [r,R]} g_{G_2} (x)
=
\begin{cases}
 {\frac{r^2}{2r^s(r + 1)},} & {s \leqslant 1} \\
 {\frac{R^2}{2R^s(R + 1)},} & {s \geqslant 2} \\
\end{cases}
\end{equation}

\noindent and
\begin{equation}
\label{eq138} M = \mathop {\sup }\limits_{x \in [r,R]} g_{G_2} (x)
=
\begin{cases}
 {\frac{R^2}{2R^s(R + 1)},} & {s \leqslant 1} \\
 {\frac{r^2}{2r^s(r + 1)},} & {s \geqslant 2} \\
\end{cases}.
\end{equation}

Now (\ref{eq137}) and (\ref{eq138}) together with (\ref{eq56})
give the required result.
\end{proof}

Some particular cases of the Theorem \ref{the46} are summarized in
the following corollary.

\begin{corollary} \label{cor46} The following bounds hold:
\begin{equation}
\label{eq139}
\frac{r^3}{4(r + 1)}\chi ^2(Q\vert \vert P) \leqslant G(Q\vert \vert P)
\leqslant \frac{R^3}{4(R + 1)}\chi ^2(Q\vert \vert P),
\end{equation}
\begin{equation}
\label{eq140}
\frac{r^2}{2(r + 1)}K(Q\vert \vert P) \leqslant G(Q\vert \vert P) \leqslant
\frac{R^2}{2(R + 1)}K(Q\vert \vert P),
\end{equation}
\begin{equation}
\label{eq141}
\frac{2r\sqrt r }{(r + 1)}h(P\vert \vert Q) \leqslant G(Q\vert \vert P)
\leqslant \frac{2R\sqrt R }{(R + 1)}h(P\vert \vert Q),
\end{equation}
\begin{equation}
\label{eq142}
\frac{r}{2(r + 1)}K(P\vert \vert Q) \leqslant G(Q\vert \vert P) \leqslant
\frac{R}{2(R + 1)}K(P\vert \vert Q)
\end{equation}

\noindent and
\begin{equation}
\label{eq143}
\frac{1}{4(R + 1)}\chi ^2(P\vert \vert Q) \leqslant G(Q\vert \vert P)
\leqslant \frac{1}{4(r + 1)}\chi ^2(P\vert \vert Q).
\end{equation}
\end{corollary}

\begin{proof} Inequalities (\ref{eq139}), (\ref{eq140}),
(\ref{eq141}) and (\ref{eq142}) follows from (\ref{eq133}) by
taking $s = - 1$, $s = 0$, $s = \frac{1}{2}$ and $s = 1$
respectively. The inequalities (\ref{eq143}) follows from
(\ref{eq134}) by taking $s = 2$.
\end{proof}

\begin{remark} \label{rem41} The inequalities (\ref{eq94}),
(\ref{eq95}), (\ref{eq112}), (\ref{eq113}), (\ref{eq128}),
(\ref{eq129}), (\ref{eq142}) and (\ref{eq143}) can be re-written
as
\begin{equation}
\label{eq144} r \leqslant \xi _t (P\vert \vert Q) \leqslant R,
\,\, t = 1,2,3,4,5,6,7\mbox{ and }8,
\end{equation}

\noindent where
\[
\xi _1 (P\vert \vert Q) = \frac{\sqrt {2F(P\vert \vert Q)}
}{\sqrt{\chi ^2(Q\vert \vert P)} - \sqrt {2F(P\vert \vert Q)} },
\]
\[
\xi _2 (P\vert \vert Q) = \frac{\sqrt {K(P\vert \vert Q)} - \sqrt
{F(P\vert \vert Q)} }{\sqrt {F(P\vert \vert Q)} },
\]
\[
\xi _3 (P\vert \vert Q) = \frac{\sqrt {F(Q\vert \vert P)} }{\sqrt {K(Q\vert
\vert P)} - \sqrt {F(Q\vert \vert P)} },
\]
\[
\xi _4 (P\vert \vert Q) = \frac{\sqrt {\chi ^2(P\vert \vert Q)} -
\sqrt {2 F(Q\vert \vert P)} }{\sqrt {2 F(Q\vert \vert P)} },
\]
\[
\xi _5 (P\vert \vert Q) = \frac{4G(P\vert \vert Q)}{\chi ^2(Q\vert
\vert P) - 4G(P\vert \vert Q)},
\]
\[
\xi _6 (P\vert \vert Q) = \frac{K(Q\vert \vert P) - 2G(P\vert \vert
Q)}{2G(P\vert \vert Q)},
\]
\[
\xi _7 (P\vert \vert Q) = \frac{2G(Q\vert \vert P)}{K(P\vert \vert
Q) - 2G(Q\vert \vert P)},
\]
\noindent and
\[
\xi _8 (P\vert \vert Q) = \frac{\chi ^2(P\vert \vert Q) - 4G(Q\vert \vert
P)}{4G(Q\vert \vert P)},
\]
respectively.
\end{remark}

\section{Bounds on Symmetric Divergence Measures}

In this section we shall obtain bound on \textit{symmetric
divergence measures} given by (\ref{eq2}), (\ref{eq10}) and
(\ref{eq11}) in terms of \textit{relative information of type s}
given by (\ref{eq14}). Some particular cases are also given.

\begin{theorem} \label{the51} The following bounds on
\textit{J-divergence} hold:
\begin{equation}
\label{eq145}
\frac{1 + r}{r^s}\Phi _s (P\vert \vert Q)
 \leqslant J(P\vert \vert Q) \leqslant
\frac{1 + R}{R^s}\Phi _s (P\vert \vert Q), \,\, s \leqslant 0
\end{equation}

\noindent and
\begin{equation}
\label{eq146}
\frac{1 + R}{R^s}\Phi _s (P\vert \vert Q)
 \leqslant J(P\vert \vert Q) \leqslant
\frac{1 + r}{r^s}\Phi _s (P\vert \vert Q), \,\, s \geqslant 1
\end{equation}
\end{theorem}

\begin{proof} Let us consider
\begin{equation}
\label{eq147} g_J (x) = x^{2 - s}{f}''_J (x) = x^{ - s} + x^{1 -
s} = \frac{1 + x}{x^s}, \,\, x \in (0,\infty ),
\end{equation}

\noindent where ${f}''_J (x)$ is as given by (\ref{eq39}).

From (\ref{eq147}) one gets
\begin{equation}
\label{eq148}
{g}'_J (x) = x^{ - s - 1}\left[ {(1 - s)x + ( - s)} \right]
\begin{cases}
 { \geqslant 0,} & {s \leqslant 0} \\
 { \leqslant 0,} & {s \geqslant 1} \\
\end{cases}.
\end{equation}

In view of (\ref{eq148}), we conclude that
\begin{equation}
\label{eq149} m = \mathop {\inf }\limits_{x \in [r,R]} g_J (x) =
\begin{cases}
 {\frac{1 + r}{r^s},} & {s \leqslant 0} \\
 {\frac{1 + R}{R^s},} & {s \geqslant 1} \\
\end{cases}
\end{equation}

\noindent and
\begin{equation}
\label{eq150} M = \mathop {\sup }\limits_{x \in [r,R]} g_J (x) =
\begin{cases}
 {\frac{R^{1 - s}}{2(1 + R)},} & {s \leqslant 0} \\
 {\frac{r^{1 - s}}{2(1 + r)},} & {s \geqslant 1} \\
\end{cases}
\end{equation}

Now (\ref{eq149}) and (\ref{eq150}) together with (\ref{eq56})
give the required result.
\end{proof}

Some particular cases of the Theorem \ref{the51} are summarized in
the following corollary.

\begin{corollary} \label{cor51} The following bounds hold:
\begin{equation}
\label{eq151}
\frac{r + r^2}{2}\chi ^2(Q\vert \vert P) \leqslant J(P\vert \vert Q)
\leqslant \frac{R + R^2}{2}\chi ^2(Q\vert \vert P),
\end{equation}
\begin{equation}
\label{eq152}
(1 + r)K(Q\vert \vert P) \leqslant J(P\vert \vert Q) \leqslant (1 +
R)K(Q\vert \vert P),
\end{equation}
\begin{equation}
\label{eq153}
\frac{1 + R}{R}K(P\vert \vert Q) \leqslant J(P\vert \vert Q) \leqslant
\frac{1 + r}{r}K(P\vert \vert Q)
\end{equation}

\noindent and
\begin{equation}
\label{eq154}
\frac{1 + R}{2R^2}\chi ^2(P\vert \vert Q) \leqslant J(P\vert \vert Q)
\leqslant \frac{1 + r}{2r^2}\chi ^2(P\vert \vert Q).
\end{equation}
\end{corollary}

\begin{proof} Inequalities (\ref{eq151}) and (\ref{eq152})
follows from (\ref{eq145}) by taking $s = - 1$ and $s = 0$
respectively. The inequalities (\ref{eq153}) and (\ref{eq154})
follows from (\ref{eq146}) by taking $s = 1$ and $s = 2$
respectively.
\end{proof}

The cases $s = \frac{1}{2}$ is not included in the inequalities
(\ref{eq145}) and (\ref{eq146}). This we shall do separately in
the following proposition.

\begin{proposition} \label{pro51} The following bound holds:
\begin{equation}
\label{eq155}
h(P\vert \vert Q)
 \leqslant \frac{1}{8}J(P\vert \vert Q).
\end{equation}
\end{proposition}

\begin{proof} Take $s = \frac{1}{2}$ in (\ref{eq147}), we
have
\begin{equation}
\label{eq156}
g_J (x) = \frac{x + 1}{\sqrt x }
\end{equation}

\noindent and
\begin{equation}
\label{eq157}
{g}'_J (x) = \frac{x - 1}{2x^{3 / 2}}
\begin{cases}
 { \geqslant 0,} & {x \geqslant 1} \\
 { \leqslant 0,} & {x \leqslant 1} \\
\end{cases}.
\end{equation}

Thus from (\ref{eq157}) we conclude that the function $g_J (x)$
given by (\ref{eq156}) is decreasing in $x \in (0,1)$ and
increasing in $x \in (1,\infty )$, and hence
\begin{equation}
\label{eq158}
m = \mathop {\inf }\limits_{x \in (0,\infty )} g_J (x) = g_J (1) = 2.
\end{equation}

Now (\ref{eq158}) together with (\ref{eq56}) give the inequality
(\ref{eq155}).
\end{proof}

\begin{theorem} \label{the52} The following bounds on
\textit{JS-divergence} hold:
\begin{equation}
\label{eq159}
\frac{r^{1 - s}}{2(1 + r)}\Phi _s (P\vert \vert Q)
 \leqslant I(P\vert \vert Q) \leqslant
\frac{R^{1 - s}}{2(1 + R)}\Phi _s (P\vert \vert Q), \,\, s
\leqslant 0
\end{equation}

\noindent and
\begin{equation}
\label{eq160}
\frac{R^{1 - s}}{2(1 + R)}\Phi _s (P\vert \vert Q)
 \leqslant I(P\vert \vert Q) \leqslant
\frac{r^{1 - s}}{2(1 + r)}\Phi _s (P\vert \vert Q), \,\, s
\geqslant 1.
\end{equation}
\end{theorem}

\begin{proof} Let us consider
\begin{equation}
\label{eq161} g_I (x) = x^{2 - s}{f}''_I (x) = \frac{x^{1 -
s}}{2(x + 1)}, \,\, x \in (0,\infty ),
\end{equation}

\noindent where ${f}''_I (x)$ is as given by (\ref{eq42}).

From (\ref{eq161}) one gets
\begin{equation}
\label{eq162}
{g}'_I (x) = \frac{x^{ - s}\left[ {( - s)x + (1 - s)} \right]}{2(x + 1)^2}
\begin{cases}
 { \geqslant 0,} & {s \leqslant 0} \\
 { \leqslant 0,} & {s \geqslant 1} \\
\end{cases}.
\end{equation}

In view of (\ref{eq162}), we conclude that
\begin{equation}
\label{eq163} m = \mathop {\inf }\limits_{x \in [r,R]} g_I (x) =
\begin{cases}
 {\frac{r^{1 - s}}{2(1 + r)},} & {s \leqslant 0} \\
 {\frac{R^{1 - s}}{2(1 + R)},} & {s \geqslant 1} \\
\end{cases}
\end{equation}

\noindent and
\begin{equation}
\label{eq164} M = \mathop {\sup }\limits_{x \in [r,R]} g_I (x) =
\begin{cases}
 {\frac{R^{1 - s}}{2(1 + R)},} & {s \leqslant 0} \\
 {\frac{r^{1 - s}}{2(1 + r)},} & {s \geqslant 1} \\
\end{cases}
\end{equation}

Now (\ref{eq163}) and (\ref{eq164}) together with (\ref{eq56})
give the required result.
\end{proof}

Some particular cases of the Theorem \ref{the52} are summarized in
the following corollary.

\begin{corollary} \label{cor52} The following bounds hold:
\begin{equation}
\label{eq165}
\frac{r^2}{4(1 + r)}\chi ^2(Q\vert \vert P) \leqslant I(P\vert \vert Q)
\leqslant \frac{R^2}{4(1 + R)}\chi ^2(Q\vert \vert P),
\end{equation}
\begin{equation}
\label{eq166}
\frac{r}{2(1 + r)}K(Q\vert \vert P) \leqslant I(P\vert \vert Q) \leqslant
\frac{R}{2(1 + R)}K(Q\vert \vert P).
\end{equation}
\begin{equation}
\label{eq167}
\frac{1}{2(1 + R)}K(P\vert \vert Q) \leqslant I(P\vert \vert Q) \leqslant
\frac{1}{2(1 + r)}K(P\vert \vert Q),
\end{equation}

\noindent and
\begin{equation}
\label{eq168}
\frac{1}{4(R + R^2)}\chi ^2(P\vert \vert Q) \leqslant I(P\vert \vert Q)
\leqslant \frac{1}{4(r + r^2)}\chi ^2(P\vert \vert Q),
\end{equation}
\end{corollary}

\begin{proof} Inequalities (\ref{eq165}) and (\ref{eq166})
follows from (\ref{eq159}) by taking $s = - 1$ and $s = 0$
respectively. The inequalities (\ref{eq167}) and (\ref{eq168})
follows from (\ref{eq160}) by taking $s = 1$ and $s = 2$
respectively.
\end{proof}

The case $s = \frac{1}{2}$ is not included in the inequalities
(\ref{eq159}) and (\ref{eq160}). This we shall do separately in
the following proposition.

\begin{proposition} \label{pro52} The following inequality hold:
\begin{equation}
\label{eq169}
I(P\vert \vert Q) \leqslant h(P\vert \vert Q).
\end{equation}
\end{proposition}

\begin{proof} For $s = \frac{1}{2}$ in (\ref{eq161}), we
have
\begin{equation}
\label{eq170}
g_I (x) = \frac{\sqrt x }{2(x + 1)}
\end{equation}

\noindent and
\begin{equation}
\label{eq171}
{g}'_I (x) = \frac{1 - x}{4\sqrt x (x + 1)^2}
\begin{cases}
 { \geqslant 0,} & {x \leqslant 1} \\
 { \leqslant 0,} & {x \geqslant 1} \\
\end{cases}.
\end{equation}

Thus from (\ref{eq171}) we conclude that the function $g_I (x)$
given by (\ref{eq170}) is increasing in $x \in (0,1)$ and
decreasing in $x \in (1,\infty )$, and hence
\begin{equation}
\label{eq172}
M = \mathop {\sup }\limits_{x \in (0,\infty )} g_I (x) = g_I (1) =
\frac{1}{4}.
\end{equation}

Now (\ref{eq172}) and (\ref{eq56}) together give the inequality
(\ref{eq169}).
\end{proof}

\begin{theorem} \label{the53} The following bounds on
\textit{AG-divergence} hold:
\begin{equation}
\label{eq173}
\frac{r^{ - s}(1 + r^2)}{4(1 + r)}\Phi _s (P\vert \vert Q)
 \leqslant T(P\vert \vert Q) \leqslant
\frac{R^{ - s}(1 + R^2)}{4(1 + R)}\Phi _s (P\vert \vert Q), \,\, s
\leqslant - 1
\end{equation}

\noindent and
\begin{equation}
\label{eq174}
\frac{R^{ - s}(1 + R^2)}{4(1 + R)}\Phi _s (P\vert \vert Q)
 \leqslant T(P\vert \vert Q) \leqslant
\frac{r^{ - s}(1 + r^2)}{4(1 + r)}\Phi _s (P\vert \vert Q), \,\, s
\geqslant 2.
\end{equation}
\end{theorem}

\begin{proof} Let us consider
\begin{equation}
\label{eq175} g_T (x) = x^{2 - s}{f}''_T (x) = \frac{x^{2 - s}(1 +
x^2)}{4(x^2 + x^3)} = \frac{x^{ - s} + x^{2 - s}}{4(1 + x)}, \,\,
x \in (0,\infty )
\end{equation}

\noindent where ${f}''_T (x)$ is as given by (\ref{eq45}).

From (\ref{eq175}) one gets
\begin{equation}
\label{eq176} {g}'_T (x) = - \frac{x^{ - s - 1}\left[ {(s - 1)x^3
+ (s - 2)x^2 + (s + 1)x + s} \right]}{4(x + 1)^2}
\begin{cases}
 { \geqslant 0,} & {s \leqslant - 1} \\
 { \leqslant 0,} & {s \geqslant 2} \\
\end{cases}
\end{equation}

In view of (\ref{eq176}), we conclude that
\begin{equation}
\label{eq177} m = \mathop {\inf }\limits_{x \in [r,R]} g_T (x) =
\begin{cases}
 {\frac{r^{ - s}(1 + r^2)}{4(1 + r)},} & {s \leqslant - 1} \\
 {\frac{R^{ - s}(1 + R^2)}{4(1 + R)},} & {s \geqslant 2} \\
\end{cases}
\end{equation}

\noindent and
\begin{equation}
\label{eq178} M = \mathop {\sup }\limits_{x \in [r,R]} g_T (x) =
\begin{cases}
 {\frac{R^{ - s}(1 + R^2)}{4(1 + R)},} & {s \leqslant - 1} \\
 {\frac{r^{ - s}(1 + r^2)}{4(1 + r)},} & {s \geqslant 2} \\
\end{cases}
\end{equation}

Now (\ref{eq177}) and (\ref{eq178}) together with (\ref{eq56})
give the required result.
\end{proof}

Some particular cases of the Theorem \ref{the53} are summarized in
the following corollary.

\begin{corollary} \label{cor53} The following bounds hold:
\begin{equation}
\label{eq179}
\frac{r + r^3}{8(1 + r)}\chi ^2(Q\vert \vert P) \leqslant T(P\vert \vert Q)
\leqslant \frac{R + R^3}{8(1 + R)}\chi ^2(Q\vert \vert P)
\end{equation}

\noindent and
\begin{equation}
\label{eq180}
\frac{1 + R^2}{8R^2(1 + R)}\chi ^2(P\vert \vert Q) \leqslant T(P\vert \vert
Q) \leqslant \frac{1 + r^2}{8r^2(1 + r)}\chi ^2(P\vert \vert Q).
\end{equation}
\end{corollary}

\begin{proof} Inequalities (\ref{eq179}) follows from
(\ref{eq173}) by taking $s = - 1$. The inequalities (\ref{eq180})
follows from (\ref{eq174}) by taking $s = 2$.
\end{proof}

The cases $s = 0$, $s = \frac{1}{2}$ and $s = 1$ are not included
in the inequalities (\ref{eq173}) and (\ref{eq174}). This we shall
do separately in the following proposition.

\begin{proposition} \label{pro53} The following inequalities hold:
\begin{equation}
\label{eq181} \frac{\sqrt 2 - 1}{2} K(Q\vert \vert P) \leqslant
T(P\vert \vert Q),
\end{equation}
\begin{equation}
\label{eq182} \frac{\sqrt 2 - 1}{2} K(P\vert \vert Q) \leqslant
T(P\vert \vert Q)
\end{equation}

\noindent and
\begin{equation}
\label{eq183}
h(P\vert \vert Q) \leqslant T(P\vert \vert Q).
\end{equation}
\end{proposition}

\begin{proof} Take $s = 0$ in (\ref{eq175}), we have
\begin{equation}
\label{eq184}
g_T (x) = \frac{x^2 + 1}{4(x + 1)}
\end{equation}

\noindent and
\begin{equation}
\label{eq185} {g}'_T (x) = \frac{(x + 1 - \sqrt 2 )(x + 1 + \sqrt
2 )}{4(x + 1)^2}\begin{cases}
 { \geqslant 0,} & {x \geqslant \sqrt 2 - 1} \\
 { \leqslant 0,} & {x \leqslant \sqrt 2 - 1} \\
\end{cases}.
\end{equation}

In view of (\ref{eq185}) we conclude that the function $g_T (x)$
given by (\ref{eq184}) is decreasing in $x \in (0,\sqrt 2 - 1)$
and increasing in $x \in (\sqrt 2 - 1,\infty )$, and hence
\begin{equation}
\label{eq186} m = \mathop {\inf }\limits_{x \in (0,\infty )} g_T
(x) = g_T (\sqrt 2 - 1) = \frac{1 + (\sqrt 2 - 1)^2}{4\sqrt 2 } =
\frac{\sqrt 2 - 1}{2}.
\end{equation}

Now (\ref{eq186}) together with (\ref{eq56}) give the inequality
(\ref{eq181}).

Again take $s = 1$ in (\ref{eq175}), we have
\begin{equation}
\label{eq187}
g_T (x) = \frac{x^2 + 1}{4x(x + 1)},
\end{equation}

\noindent and
\begin{equation}
\label{eq188} {g}'_T (x) = \frac{(x - 1 - \sqrt 2 )(x - 1 + \sqrt
2 )}{4(x + 1)^2}\begin{cases}
 { \geqslant 0,} & {x \geqslant \sqrt 2 + 1} \\
 { \leqslant 0,} & {x \leqslant \sqrt 2 + 1} \\
\end{cases}.
\end{equation}

In view of (\ref{eq188}), we conclude that the function $g_T (x)$
given by (\ref{eq187}) is decreasing in $x \in (0,\sqrt 2 + 1)$
and increasing in $x \in (\sqrt 2 + 1,\infty )$, and hence
\begin{equation}
\label{eq189} m = \mathop {\inf }\limits_{x \in (0,\infty )} g_T
(x) = g_T (\sqrt 2 + 1) = \frac{\sqrt 2 - 1}{2}.
\end{equation}

Now (\ref{eq189}) together with (\ref{eq56}) give the inequality
(\ref{eq182}).

Finally, take $s = \frac{1}{2}$ in (\ref{eq175}), we have
\begin{equation}
\label{eq190}
g_T (x) = \frac{x^2 + 1}{4\sqrt x (x + 1)},
\end{equation}

\noindent and
\begin{equation}
\label{eq191}
{g}'_T (x) = \frac{(x - 1)\left[ {(x - 1)^2 + 6x} \right]}{8x^{3 / 2}(x +
1)^2}
\begin{cases}
 { \geqslant 0,} & {x \geqslant 1} \\
 { \leqslant 0,} & {x \leqslant 1} \\
\end{cases}.
\end{equation}

In view of (\ref{eq191}) we conclude that the function $g_T (x)$
given by (\ref{eq190}) is decreasing in $x \in (0,1)$ and
increasing in $x \in (1,\infty )$, and hence
\begin{equation}
\label{eq192}
m = \mathop {\inf }\limits_{x \in (0,\infty )} g_T (x) = g_T (1) =
\frac{1}{4}.
\end{equation}

Now (\ref{eq192}) together with (\ref{eq56}) give the inequalities
(\ref{eq183}).
\end{proof}

\begin{remark} \label{rem51}
\begin{itemize}
\item[(i)] The inequalities appearing in (\ref{eq152}),
(\ref{eq153}), (\ref{eq166}) and (\ref{eq167}) can be re-written
as
\begin{equation}
\label{eq193} r \leqslant \zeta _t (P\vert \vert Q) \leqslant R,
\,\, t = 1,2,3\mbox{ and 4},
\end{equation}

\noindent where
\[
\zeta _1 (P\vert \vert Q) = \frac{J(P\vert \vert Q) - K(Q\vert \vert
P)}{K(Q\vert \vert P)},
\]
\[
\zeta _2 (P\vert \vert Q) = \frac{K(P\vert \vert Q)}{J(P\vert \vert Q) -
K(P\vert \vert Q)},
\]
\[
\zeta _3 (P\vert \vert Q) = \frac{2I(P\vert \vert Q)}{K(Q\vert \vert P) -
2I(P\vert \vert Q)},
\]

\noindent and
\[
\zeta _4 (P\vert \vert Q) = \frac{K(P\vert \vert Q) - 2I(P\vert \vert
Q)}{2I(P\vert \vert Q)}
\]

\noindent respectively.

\item[(ii)] In view of (\ref{eq155}), (\ref{eq169}) and
(\ref{eq183}), we have the following interesting relation
\begin{equation}
\label{eq194} I(P\vert \vert Q) \leqslant h(P\vert \vert Q)
\leqslant \left( {T(P\vert \vert Q)\mbox{ or }\frac{1}{8}J(P\vert
\vert Q)} \right)
\end{equation}

A general form of the inequalities (\ref{eq194}) can be seen in
Taneja \cite{tan6} , where more kind of symmetric measures are
also studied.
\end{itemize}
\end{remark}


\begin{thebibliography}{99}
\bibitem{pea} K. PEARSON, On the Criterion that a given System of Deviations
from the Probable in the case of Correlated System of Variables is
such that it can be Reasonable supposed to have arisen from Random
Sampling, \textit{Phil. Mag.}, \textbf{50}(1900), 157-172.
\end{thebibliography}
\end{document}